\input amstex
\magnification 1200
\documentstyle{amsppt}
\pageheight {7.291in}
\pagewidth {5.5in}
\NoRunningHeads
\NoBlackBoxes

\def\barr{\begin{displaymath} \begin{array}{lll}}
\def\earr{\end{array} \end{displaymath}}
\def\be{\begin{eqnarray}}
\def\ee{\end{eqnarray}}
\def\n{\nonumber }
\def\oo{\otimes}
\def\la{\lambda}
\def\al{\alpha}
\def\bt{\beta}
\def\dt{\delta}
\def\la{\lambda}

\def\si{\sigma}
\def\gm{\gamma}
\def\td{\tilde \delta}
\def\h{\eta}

\def\tl{\triangleleft}
\def\tr{\triangleright}

\rightheadtext{Associative triples}
\leftheadtext{A.I.Mudrov}
\email {aimudrov\@dg2062.spb.edu}
\endemail
\thanks
I thank A. A. Stolin who drew my attention 
to a link between Yang-Baxter equation and
a cyclic symmetric inner product. I am grateful to 
P. P. Kulish for discussions.
\endthanks

\document

\def\G{\hat \Gamma}
\def\B{{\overline B}}
\def\ell{{\text{ell}}}
\def\Tr{\text{Tr}}
\def\Mat{\text{Mat}}

\def\id{\text{id}}
\def\card{\text{card}}

\def\M{\frak M}
\def\N{\frak N}
\def\D{\frak D}
\def\A{\frak A}
\def\B{\frak B}

\def\P{\Cal P}

\def\be{\bold e}

\def\R{\frak R}
\def\h1{\hat{\bold 1}}

\def\a{\frak a}

\def\Ua{U_q(\tilde\g)}
\def\U2{{\Ua}_2}
\def\g{\frak g}
\def\n{\frak n}

\def\C{\Bbb C}
\def\d{\frak d}

\def\Tr{\text{\rm Tr}}
\def\l{\frak l}

\def\<{\langle}
\def\>{\rangle}

\def\op{\oplus}

\def\h{{\frak h}}
\baselineskip 12 pt

\topmatter
\title
Associative triples and Yang-Baxter equation
\endtitle
\author {\rm {\bf Andrei Mudrov}}\\
Department of Mathematics, Bar-Ilan University, 52900 Ramat-Gan,
Israel;\\
e-mail: mudrova\@macs.biu.ac.il
\endauthor
\endtopmatter
\noindent
$$\vbox{\hsize 23 pc\noindent
\strut
{
\eightrm
\vskip .03in
\noindent
ABSTRACT. We introduce triples of associative algebras as a tool for
building solutions to the Yang-Baxter equation.
It turns out that the class of R-matrices thus obtained
is related to a Hecke-like condition, which is  formulated for associative
algebras with symmetric cyclic inner product.
R-matrices for a subclass of the $A_n$-type Belavin-Drinfel'd triples
are derived in this way.\newline
\vskip 12pt
\noindent
AMS classification codes: 81R50, 16W99.
}
\strut}
$$

\vskip .3in

\subhead 1. Introduction
\endsubhead

The canonical Faddeev-Reshetikhin-Takhtajan recipe for constructing
quantum matrix groups  \cite{RTF} is based on a solution
to the Yang-Baxter equation (YBE), which is assumed to be {\it a priori}
known.
On the other hand, the theory of quantum groups  was
designed as an environment for constructing such solutions \cite{D}
being of interest for mathematics and physics.
The quasi-classical limit of YBE, the classical Yang-Baxter
 equation (cYBE), has clear algebraic interpretation in
terms of Manin triples, which have been classified for
semisimple Lie algebras in \cite{BD}.
Transition from  Manin triples to quantum universal
enveloping algebras brings about the notions
of quantum double and universal solution  of YBE
\cite{D}.
The possibility to quantize an arbitrary Manin triple has
been proven in \cite{EK}, and for universal R-matrices 
there were derived rather complicated
although explicit formulas, in the semisimple case, \cite{ESS}.  
At the same time,
one can raise the question what structure particulars of an
associative ring itself may be responsible for YBE,
without appealing to
the intricate technique of the Hopf algebra
theory. 
An explicit formula for $sl_n(\C)$ R-matrices, known 
as the GGS conjecture, was proposed in \cite{GGS}. It
 has been confirmed \cite{GH,Sch1} in
many cases and a combinatorial proof was
recently found in \cite{Sch2}.
In the present paper, we pursue another
approach making use of cyclic inner product on associative
algebras and
keeping an analogy with Manin triples. Such
a point of view has led us to definition of
associative triples, which allow to explain many
examples and quantize a wide class of Belavin-Drinfel'd triples
associated to the special linear Lie algebras.
In associative triples,  
R-matrices  naturally  naturally  
into  the sum of two solutions to YBE, one of
them being a part of the canonical element of
the cyclic inner product and the other
belonging to a smaller subalgebra.
That "smaller" solution should satisfy a Hecke-like
condition.
The problem of building R-matrices in a given algebra
is thus reduced to finding an associative triple (if that possible)
with its total algebra as the homomorphic pre-image.
For $\Mat_n(\C)$, we take
a proper extension of $\Mat_n(\C)\op \Mat_n(\C)$
for the role of such a pre-image.

\subhead 2. Cyclic inner product and YBE
\endsubhead
Throughout the paper, we assume ${\M}$ to be a finite
dimensional  associative algebra with unit over $\C$.
Nevertheless, there are interesting infinite-dimensional
examples that can be understood in the framework of our 
construction. 
Arbitrary fields are also admissible, and presence of unit 
is not necessary in certain cases.

The Yang-Baxter equation in $\M^{\oo 3}$ is
$$
R_{12}R_{13}R_{23}=R_{23}R_{13}R_{12},
\tag 1
$$
where $R\in \M^{\oo2}$ and the subscripts specify
the way of embedding $\M^{\oo 2}$ into  $\M^{\oo 3}$.
Supposing R being a deformation of the unit, 
$R=1\oo 1 +\la r + o(\la)$, $\la \in \C$,
the element $r$ satisfies the classical Yang-Baxter equation
$$
[r_{12},r_{13}]+  [r_{12},r_{23}]+  [r_{13},r_{23}]=0.
\tag 2
$$
The square brackets mean the commutator
$[\al_1,\al_2]=\al_1 \al_2-\al_2 \al_1$, $\al_i \in \M$.

Suppose  $\M$ is endowed  with a non-degenerate symmetric
cyclic inner product $(\cdot,\cdot)$, that is
$(\al_1\al_2,\al_3) = (\al_2,\al_3\al_1)$ for all $\al_i \in \M$.
Choose a basis $\{\al_i\}\in \M$ and let $\{\al^i\}$ be its dual:
$(\alpha_i , \al^k)=\delta_i^k$ (the Kronekker symbols).

\vskip .03in
\demo{Definition 1}
The element $\si_{\M} =\al_i \oo \al^i\in \M^{\oo 2}$
(summation over repeating upper and lower indices understood 
throughout the paper)
is called {\it permutation} in $\M$.
\enddemo
\noindent
For unital algebras, cyclic inner products are in one-to-one 
correspondence with linear  functionals $t_\M$ obeying 
$t_\M(\al\bt)=t_\M(\bt\al)=(\al,\bt)$.
In  the  endomorphism algebra of the
vector space $\C^n$ with the trace functional, the permutation is
the flip operator: $x\oo y\to y\oo x \in \C^n\oo \C^n$.
\proclaim{Proposition 1 {\rm(}\cite{BFS}{\rm )}}
Permutation satisfies YBE.
For all $\al\in\M$
$$(\al\oo 1)\si_{\M} = \si_{\M}(1\oo \al),
\quad (1\oo \al)\si_{\M} = \si_{\M} (\al\oo 1).$$
\endproclaim
\vskip .03in
\noindent
Remark that associative algebras with a functional $t$ such
that the  bilinear form $\M\oo \M \to k$, 
$\al_1\oo\al_2 \to t(\al_1\al_2)$, is non-degenerate are called
Frobenius. For those algebras, the canonical element fulfills
the braid relation too, in analogy with the conventional 
matrix permutation \cite{BFS}.

Recall that a Manin triple $(\g,\a,\a^*)$
comprises a Lie algebra $\g$  with an ad-invariant
inner product and its two lagrangian (maximal isotropic) Lie subalgebras $\a$
and $\a^*$ with zero intersection.
The canonical element $\alpha_i \otimes \beta^i \in \a\otimes\a^*$
solves cYBE.
To find a possible "quantification"  of this construction,
 consider first the situation when a Manin triple
$(\M,\A,\A^*)$ is formed by the commutator Lie algebras of associative
algebras, and the ad-invariant 2-form
is induced by a non-degenerate cyclic inner product.
Consider $\A$ and $\A^*$ as bimodules for each other, the
left $(\tr)$ and right  $(\tl)$ actions being
induced via duality from 
multiplication.
The  product in $\M$
is expressed by the formula
$$
\al \bt = \al \tr \bt + \al \tl \bt , \quad
\bt \al = \bt \tr \al + \bt \tl \al , \quad
\al \in \A
, \quad
\bt\in \A^*.
\tag 3
$$
Associativity is encoded in the following two equations:
$$
  (\al_1\tr \bt_1, \al_2 \tl \bt_2) + ( \al_1\tl \bt_1 , \al_2 \tr
\bt_2)
= (\bt_1\tr \al_2, \bt_2 \tl \al_1) + ( \bt_1\tl \al_2 , \bt_2 \tr
\al_1)
,
$$
$$
 ( \al_2 \al_1,\bt_1 \bt_2 )=
 ( \bt_2 \tr \al_2, \al_1\tr\bt_1) +
 ( \al_1 \tl \bt_1,\bt_2 \tl \al_2 )
, \quad \al_i \in \A, \quad \bt_i \in \A^*.
$$
Conversely, the linear sum $\A\oplus\A^*$ is equipped
with the structure of an associative algebra (3), 
provided the actions $\tl$ and $\tr$  satisfy
these two conditions. Then, the following statement holds true.
\proclaim{Proposition 2}
Natural pairing between $\A$ and $\A^*$ induces a non-degenerate
cyclic inner product on $\M=\A\oplus\A^*$ such that $\A$ and $ \A^*$
are isotropic.
The canonical element $\al_i\oo \al^i\in \A\oo \A^*\subset \M^{\oo 2}$
of the pairing satisfies YBE.
\endproclaim
\vskip .03in
\demo{Proof}
The first statement is immediate.
Verification of the second is conducted in the next section 
in a more general setting
$\blacksquare$\enddemo

For the above construction,  unit is not necessary. But if it
is present, the sum 
$$
1\oo 1 + \la \al_i \oo \al^i
\tag 4
$$
is also a solution to YBE for an arbitrary value of the scalar parameter
$\la$ (cf. concluding remarks).
This is true because  $(\M,\A,\A^*)$ is a Manin triple of the 
corresponding commutator Lie algebras.
Associativity imposes strong restrictions on the algebras
$\A$ and $\A^*$; for instance, they cannot be simultaneously unital
or they would consist of the zero elements only.
Nevertheless, there is a non-trivial example, rather infinite dimensional, 
describing the XXX-spin chains in the theory of  integrable models 
\cite{F,KS}. That is  the Yang matrix
$$
R(z,u)=1\oo1 + \la \frac{ \P}{z-u} ,
$$
where $\P$ is the conventional permutation operator acting on
$\C^n\oo \C^n$.
The element $\frac{ 1}{z-u}$ is a brief way of writing the formal power
series
 $\sum_{k\geq 0}\frac{u^k}{z^{k+1}}$. It  represents
 the canonical element of pairing Res$_0$ between
 $\C[z]$ and $\frac{1}{z}\C[\frac{1}{z}]$.
In this example, the algebra $\M$ is formed by the Laurent polynomials with
matrix coefficients, but $R(z,u)$ requires 
extension by the Laurent series. 

Propositions 1 and 2 can be understood with the help of the following 
construction.
Let $\N_\pm$ be two linear subspaces in $\M$ paired
via the inner product. There exist two bijections
$\N_\pm \to \N^*_\mp$ inducing projectors 
$\pi_\pm\colon\M \to \N_\pm$. If $\{\al_i\}$ is a basis 
in $\N_+$ and $\{\bt^i\}$ its dual in $\N_-$,
then $\pi_+$ is given by $\mu\to \al_i (\mu,\bt^i)$.
The projector 
$\pi_-$, the conjugate to $\pi_+$ with respect to the 
inner product, acts as  $\mu\to (\al_i,\mu)\bt^i$.
Denote 
$\M_\pm=\{\mu|\pi_\pm(\al_1 \mu) \al_2 =\al_1\pi_\pm(\mu \al_2),  \al_1,\al_2 \in \N_\pm\}$.
 The subspaces $\M_\pm$
 contain, in particular, the normalizers for $\N_\pm$, i.e. the 
maximal subalgebras in $\M$ for which $\N_\pm$ are bimodules.
In particular, $\N_\pm\subset \M_\pm$ if $\N_\pm$ are subalgebras.
\proclaim{Proposition 3}
Suppose $\M=\M_+ + \M_-$ as a linear space. The canonical element 
$Q=\al_i\oo \bt^i\in \N_+\oo \N_- \subset \M^{\oo2}$ satisfies YBE.
\endproclaim
\vskip .03in
\demo{Proof}
 It will be sufficient to  evaluate the middle tensor components 
of YBE   on
elements of $\M_+$ and $\M_-$ separately, since they span entire $\M$.
For arbitrary $\mu\in\M_+$, the Yang-Baxter equation
$$
\al_i \al_j \oo (\bt^i \al_k,\mu) \oo \bt^j \bt^k=
\al_j \al_i \oo (\al_k \bt^i,\mu) \oo \bt^k \bt^j
$$
is rewritten as 
$$
\al_i \al_j \oo (\bt^i,\al_k \mu) \oo \bt^j \bt^k=
\al_j \al_i \oo (\bt^i,\mu \al_k) \oo \bt^k \bt^j
$$
and therefore
$$
\pi_+(\al_k \mu) \al_j \oo \bt^j \bt^k=
\al_j \pi_+ (\mu \al_k)\oo \bt^k \bt^j.
$$
Similarly one verifies YBE by pairing with
elements of $\M_-$
$\blacksquare$\enddemo

The above statement is an apparent generalization of Propositions 1 and 2.
In the first case, one has $\N_\pm=\M_\pm = \M$ and, in the second,
already $\N_+ +\N_-$ gives entire $\M$. 
Let us present another  example, quite exotic, where 
the normalizers
for $\N_\pm$ are very small whereas both $\M_\pm$ coincide with
the whole $\M$.
Take $\M$ to be the algebra of $n\times n$ matrices  
and denote  $e^i_k$  the  matrix units. 
Let $\si$ be a permutation of the set of indices $I=\{1,...,n\}$. 
Put $\N_+=\text{Span}\{ e_i^{\si(i)}\}_{i\in I}$
and 
$\N_-=\text{Span}\{e^i_{\si(i)}\}_{i\in I}$;
then  the canonical element with respect to the trace pairing reads 
$Q= \sum_{i\in I} e_i^{\si(i)}\oo e^i_{\si(i)}$. 
Now observe that $\M_+=\M$.
Indeed, for any matrix $u=u^k_i e^i_k\in \M$ 
one has $\pi_+\Bigl(u e_i^{\si(i)}\Bigr)=u^{\si(i)}_{\si(i)} e_i^{\si(i)}$
and     $\pi_+\Bigl(e_i^{\si(i)} u\Bigr)=u^i_i e_i^{\si(i)}$.
So one  gets the identity
$$
\pi_+\Bigl(e_i^{\si(i)}u\Bigr)e_k^{\si(k)}=
u^i_i \delta_i^{\si(k)} 
e_k^{\si^2(k)}\quad= 
\quad u^{\si(k)}_{\si(k)} \delta_i^{\si(k)} e_k^{\si^2(k)} =
e_i^{\si(i)}\pi_+\Bigl(u e_k^{\si(k)}\Bigl)
.
$$
Note that the normalizers for $\N^\pm$ consist of diagonal
matrices only.

Proposition 3 supplies one  with solutions that seem to 
be quite distant from those related to quantum groups.
We view it  as 
a tool for constructing parts of R-matrices of
real interest.
\subhead 3. Associative triples
\endsubhead

From now on we assume the  subspaces $\N_\pm$ to be subalgebras.
Suppose  $\M$ can be represented as the linear
sum of three subalgebras
$\N_- + \D + \N_+$ such that \newline 
{ i)} $\N_-$ is dual to $\N_+$ with 
respect to the inner product,\newline
{ ii)} $\N_\pm$ are $\D$-bimodules,\newline
{ iii)} $\D$ is orthogonal to $\N_- +\N_+$.\newline
In the sequel, $\M_\pm $ will stand for the subalgebras  $\D + \N_\pm$.
\vskip .03in
\demo{Definition 2}
The set $(\M,\M_+\M_-)$ is called {\it associative triple}
with {\it diagonal} $\D=\M_+\cap \M_-$.
A triple is called {\it trivial} if it coincides with
its diagonal and {\it disjoint} if $\D=\{0\}$.
\enddemo
\noindent
It follows from condition iii) that restriction of
the cyclic form to $\D$ is nondegenerate.
Product of two elements
from $\N_\pm\N_\mp$ may contribute to every part
of $\M$, in analogy with (3), including the diagonal. 
Thus associative triples generalize decomposition of
$\Mat_{n}(\C)$ into the diagonal, strictly upper and lower
triangular matrix subalgebras.
Let us give one more definition before
formulating the  basic statement of the paper.
\vskip .03in
\demo{Definition 3}
An element $S\in \M^{\oo2}$, where  
$\M$ is an algebra with a non-degenerate
inner product, is said to satisfy 
{\it Hecke condition} in $\M^{\oo 2}$ if 
$$
 S_{21} S-\si_{\M} S =\la^2 
\tag 5
$$
with some scalar $\la$. 
\enddemo
It is convenient
to represent  $\frac{1}{\la} = \omega=q-q^{-1}$ assuming
$q^2\not=1$ and $\la\not = 0$.
\noindent
For a  matrix algebra $\M$, this is the 
conventional Hecke condition. Moreover, that is 
the only case when the permutation
$\si_{\M}$ is invertible  \cite{BFS}. Then,
$\si^2_\M=1$ and one can combine $S$ with
$\si_{\M}$ getting a close quadratic 
equation on $\si_\M S$.
\proclaim{Theorem 1}
The canonical element $Q\in \N_+ \oo \N_-$ with
respect to the cyclic inner product satisfies YBE.
The element $R=S+Q\in \M^{\oo 2}$, where $S\in \D^{\oo 2}$,
is a solution to  YBE  if
$S$ is a solution to YBE fulfilling
the Hecke condition in $\D^{\oo 2}$.
\endproclaim
\vskip .03in
\demo{Proof}
The first assertion follows from Proposition 3 because
$\N_\pm$ are bimodules for $\M_\pm$ and $\M_+ +\M_-=\M$.
Verification of YBE is reduced to checking
$$
S_{12}Q_{13}Q_{23}+Q_{12}S_{13}Q_{23}+Q_{12}Q_{13}S_{23}+S_{12}Q_{13}S_{23}
$$
$$
=
$$
$$
Q_{23}Q_{13}S_{12}+Q_{23}S_{13}Q_{12}+S_{23}Q_{13}Q_{12}+S_{23}Q_{13}S_{12}.
$$
Indeed,  $S$ and $Q$ themselves are solutions to 
YBE and $Q$ interacts with elements of $\D$ like the permutation:
$(\dt_1\oo \dt_2)Q=Q(\dt_2\oo \dt_1)$  for every $\dt_1, \dt_2\in \D$. 
Therefore, equalities
of the type $SSQ=QSS$ hold identically and thus drop from  
YBE for the sum $S+Q$. Put  $S=\dt_i\oo \td^i \in \D^{\oo 2}$
and rewrite the equation above as
$$
\dt_i \al_j \oo \td^i \al_k \oo \bt^j \bt^k + 
\al_i \dt_j \oo \bt^i \al_k \oo \td^j \bt^k +
\al_i \al_j \oo \bt^i \dt_k \oo \bt^j \td^k +
\dt_i \al_j \oo \td^i \dt_k \oo \bt^j \td^k 
$$
$$
=
$$
$$
\al_j \dt_i \oo \al_k \td^i \oo \bt^k \bt^j +
\dt_j \al_i \oo \al_k \bt_i \oo \bt^k \td^j +
\al_j \al_i \oo \dt_k \bt_i \oo \td^k \bt^j +
\al_j \dt_i \oo \dt_k \td^i \oo \td^k \bt^j \quad
$$
Since $\M$ is spanned by $\N_\pm$ and $\D$, it suffices 
to check this identity separately pairing the middle
tensor component with their elements.

\noindent
Step 1. Pairing with $\al\in \N_+$.
$$
0 + \al_k \al \dt_j \oo \td^j \bt^k + \dt_k \al \al_j \oo \bt^j \td^k + 0 =
0 + \dt_j \al \al_k \oo \bt^k \td^j + \al_j \al \dt_k \oo \td^k \bt^j + 0.\quad
$$
Step 2. Pairing with $\bt\in \N_-$.
$$
\dt_i \al_j \oo \bt^j \bt \td^i + \al_i \dt_j \oo \bt^j \bt \td^i + 0 + 0 =
\al_j \dt_i \oo \td^i \bt \bt^j + \dt_j \al_i \oo \bt^i \bt \td^j + 0 + 0.\quad
$$
Step 3. Pairing with $\dt\in \D$.
The first and third terms on each sides turn zero.
In the fourth terms, we perform the substitution $S\to \si_\D$
of the last factors employing  
the Hecke condition. For example, 
$$
S_{12}Q_{13}S_{23}=S_{12}S_{21}Q_{13}=(S_{12}(\si_\D)_{12}+\la^2) Q_{13}=
S_{12}Q_{13}(\si_\D)_{23}+\la^2 Q_{13}.
$$
The last term will appear on both sides of the equation and
thus vanish.
The resulting equation is  
$$
0 + \al_i \dt_j \oo \td^j \dt \bt^i + 0 + \dt_i \al_j \oo \bt^j \dt \td^i =
0 + \dt_j \al_i \oo \bt^i \dt \td^j + 0 + \al_j \dt_i \oo \td^i \dt \bt^j,
$$
and it holds identically
$\blacksquare$\enddemo
\noindent
Note that representation of the  Cremmer-Gervais R-matrix
\cite{CG}
by the sum of two solutions to YBE was used in \cite{H1}.
Introduced in this paper, associative triples is an algebraic
scheme adopted to solving  the specific system of equations
arising from such a representation.

Associative triples form a category ${\Cal{AT}}$, with the 
subcategory ${\Cal{AT}}_0$ of trivial triples.
Morphisms in ${\Cal{AT}}$
are algebraic maps preserving elements of triples.
The category ${\Cal{AT}}$ admits the following operations with
its objects:
\demo{1. Transposition $\quad\M'$}
$$
(\M,\M_+,\M_-)'=(\M,\M_-,\M_+)
,\quad
t_{\M'}=t_{\M}.
$$
\enddemo
\demo{2. Sum $\quad\M^1\op \M^2$}
$$ (\M^1\op\M^2,\M^1_+\op\M^2_+,\M^1_-\op\M^2_-)
,\quad
t_{\M^1\op \M^2}=
t_{\M^1}\op t_{\M^2}.
$$
\enddemo
\demo{3. Product by objects of ${\Cal{AT}}_0 $}
$$(\A\oo\M,\A\oo\M_+,\A\oo\M_-)
,\quad
t_{\A\oo\M}=t_{\A}\oo t_{\M}.
$$
\enddemo
\demo{4. Double  $\quad (D(\M), D(\M)_+,D(\M)_-$)}
$$ D(\M)=\D\op\M\op \M,
$$
$$ D(\M)_+=\D\op\M\op\id(\M),\quad
   D(\M)_-=\D\op\D\op\id(\D)+\{0\}\op\N_+\op \N_-,
$$
$$
t_{D(\M)}=t_{\D}\op t_{\M}\ominus t_{\M}.
$$
In the definition of $D(\M)_-$ we identify
the first and last copies of $\D$.
The diagonal here is $\D\op\id(\D)\op\id(\D)$.
The functional $t_\D$ is induced by 
 $t_\M$,  and the restriction of $t_{D(\M)}$
to the last addend coincide with $-t_\M$ (that is
reflected by the notation).
\enddemo
\demo{5. Skew double  $\quad S(\M)$}
This is a disjoint triple $\M+\M^*$ of an algebra
and its dual equipped with nil multiplication.
Formula (3) degenerates into the dual regular
actions of $\M$ on $\M^*$.
\enddemo
Introduction of associative triples is motivated by the 
idea of solving  YBE in  a smaller algebra and to 
create a tool for the induction process.
Propositions 1 and 2 describe two extreme
cases of trivial and disjoint triples
providing quite exotic examples.
To find more interesting solutions, one has 
to admit non-trivial
 $\D$ and $\N_\pm$ simultaneously. 
Such applications are considered in the rest of the paper,
and this section is finished with the following statement.
\proclaim{Proposition 4}
Let $\M$ be an associative triple with the diagonal 
$\D$ and $S$ satisfies the Hecke condition in 
$\D^{\oo 2}$. Then, the element $R=S+Q$ satisfies the Hecke condition in 
$\M$ if and only if $\si_\D Q +Q^2=0$.
\endproclaim
\demo{Proof}
Taking into account $S_{12} Q=QS$ and $\si_\M=\si_\D+Q+Q_{21}$,
one has 
$$
R_{21}R=S_{21}S+S_{21}Q +QS+Q_{21}Q=\la(1\oo 1)+ \si_\M S+Q_{21}Q=
\la(1\oo 1)+ \si_\M R -\si_\D Q- Q^2
$$
as required.
$\blacksquare$\enddemo
\subhead 4.
Examples
\endsubhead
\vskip .03in
\demo{1. Drinfeld-Jimbo R-matrix for $sl_n(\C)$  \cite{J}}
\enddemo
\noindent
One can build, by recursion,
$$R_n=q\sum_{i=1}^n e^i_i\oo e^i_i +
    \sum_{i,j=1\atop i\not=j}^n e^i_i\oo e^j_j +
    \omega\sum_{i<j=1}^n e^i_j\oo e^j_i, \quad
    \omega =q-q^{-1},
$$
if assumes $\M=\Mat_{n}(\C)$ with the ordinary matrix trace,
$\M_+=\sum_{i,j=1}^{n-1}\C e^i_j +\sum_{i=1}^n\C e^i_n$,
and
$\M_-=\sum_{i,j=1}^{n-1}\C e^i_j +\sum_{i=1}^n\C e^n_i$.
So the subalgebra $\N_+$ is formed by the last matrix column
without the bottom entry. Its dual
$\N_-$ is spanned by the bottom matrix line except
the right-most diagonal element. The diagonal subalgebra $\D$
is $\Mat_{n-1}(\C)\oplus \C e_n^n$. The one-dimensional 
R-matrix $R_1=q e^n_n\oo e^n_n$ fulfills the Hecke condition in 
$\Mat_1(\C)$. By induction assumption, that holds
for the matrix $R_{n-1}$.
The direct sum $R_{n-1}+R_1$ satisfies the Yang-Baxter equation
but not the Hecke condition since the unit matrix 
$1_n\oo 1_n$ is not equal to the sum $1_{n-1}\oo 1_{n-1}+1_1\oo 1_1$. 
To fix the situation, we
put $\omega S=R_{n-1}+R_1+P_{n-1}\oo P_1+P_1\oo P_{n-1}$,
where  $P_{n-1}$ and $P_1$ stand for the projectors to the corresponding 
matrix subalgebras. Thus defined, $S$ solves the Yang-Baxter equation
too.  The  matrix  $Q=\sum_{i=1}^{n-1}e^i_n\oo   e^n_i$   fulfills   the
condition of Proposition 4, so $S+Q=\frac{1}{\omega} R_n$ is the Hecke 
matrix by induction. 

Let us pursue another represention of 
$\Mat_{n}(\C)$ as an associative triple,
taking
$\M_\pm$ to be the  subalgebras of upper and lower triangular
matrices. Then $\D$ is a commutative algebra formed by the diagonal matrices.
The Hecke condition  on
$S=a^{ik} e^i_i\oo e^k_k  \in\D^{\oo2}$, which
is an apparent solution to YBE,
boils down to the system
$$
a^{ii}+a^{ik}a^{ki}= a^{ii}a^{ii}, \quad k\not = i,\quad
a^{ik}a^{ki} = a^{im}a^{mi}, \quad k\not = i, \quad m\not = i.
\tag 6
$$
It is satisfied by the numbers $a^{ii} =\pm q^{\pm 1}/\omega$
and  $a^{ik} = b^{ik}/ \omega $,
$b^{ik}b^{ki}=1$ for $i\not = k$ (note that the classical limit 
$R\to 1$
exists, after the proper rescaling by $\omega$, only in case of  
$a^{ii} =  q/\omega$).
This is the principal solution of the Hecke condition in the
commutative
algebra $\C^n\oo\C^n$. It deviates from the standard solution $b^{ik}=1$ exactly
by the Reshetikhin diagonal  twist \cite{R}. In general,
let $\tilde S = F S F^{-1}_{21}$, where $F$ is an invertible
element from $\D^{\oo2}$, be a solution to YBE too. Then, the
Hecke condition is fulfilled and $F Q F^{-1}_{21}=Q$.
So $\tilde R = F (S+Q) F^{-1}_{21}=\tilde S+Q$ is again an
R-matrix for $\M$.
\vskip .05in
\demo{2. Baxterization procedure}
\enddemo
\noindent
The baxterization operation converts a
constant matrix solution $R$ to YBE to
that with spectral parameter
$$
R(z,u)=zR - u R^{-1}_{21},
\tag 7
$$
provided $R$ satisfies the conventional Hecke condition
$$
(R\P)^2=\omega (R\P) + I
$$
with the matrix permutation $\P$ .
Parameters $z$ and $u$ are usually represented in
the exponential form; then (7) is the
trigonometric solution to YBE. Set $\N_+=z\Mat_{n}(\C)[z]$,
$\N_-=\frac{1}{z}\Mat_{n}(\C)[\frac{1}{z}]$, and
$\D=\Mat_{n}(\C)$.
The cyclic inner product on   $\M$ is given by
the formula $(A z^k, B z^m)=\Tr(AB) \delta^{k,-m}$.
Thus we built an associative triple on $\Mat_{n}(\C)[z,\frac{1}{z}]$. Now, put 
$S=\frac{R}{\omega}$ and $Q=\frac{u}{z-u}{\P}$. The result will be
proportional to (7) because $\omega \P = R - R^{-1}_{21}$.
Again, as with the Yang matrix, one has to extend 
the Laurent polynomial algebra to that of the  Laurent series.

\subhead 5.
On quantization of Belavin-Drinfel'd triples for $sl_{n-1}(\C)$
\endsubhead

Consider  a semisimple Lie algebra $\g$ with the Cartan
subalgebra $\h$ and the polarization $\g=\n_-+\h+\n_+$.
Let $\Delta$ and $\Delta^\pm$ be, respectively,
the systems of its all, positive, and negative roots.
Recall \cite{BD} that non-skew-symmetric classical r-matrices
associated to $\g$
are  in one-to-one correspondence with Belavin-Drinfel'd  (BD)
triples
$(\Gamma_1,\Gamma_2,\tau)$ consisting of two subsets
of positive simple roots $\Gamma_i$, $i=1,2$,
and a bijection $\tau \colon\Gamma_1 \to \Gamma_2$
preserving lengths of the roots with respect to  the Killing form. 
Besides,
for every $\al\in \Gamma_1$ there is a positive integer $k$
such that $\tau^{k-1}(\al)\in \Gamma_1$ and
$\tau^{k}(\al)\not \in \Gamma_1$.
It follows from \cite{RS} (see also \cite{BZ}) that for any BD triple
the  double Lie algebra $D(\g)$ is isomorphic to the direct
sum $\g\oplus \g$ with  the invariant scalar product
being the difference of Killing forms on the two addends. 
There is a geometric description of $\g^*$ as a Lie
subalgebra in $D(\g)$  \cite{S}.
Let $\Delta_i$ be  the subsystems of roots generated by $\Gamma_i$,
$i=1,2$.
Denote $\g_i$ the semisimple Lie algebras with
the root systems  $\Delta_i$, and the 
Cartan subalgebras  $\h_i$. The bijection $\tau$ is extended
to a Lie algebra isomorphism such that
$\g_1\op \tau(\g_1)$ has trivial intersection
with the diagonal embedding $\g\op\id(\g)$. Choose  subspaces
$\h^0_i\subset (\h_i)^\perp \subset \h$ containing
their orthogonal complements in  $(\h_i)^\perp$ such that the
equation $\eta_1+\eta^0_1=\tau(\eta_1)+\eta^0_2$,
$\eta_1\in \h_1$, $\eta^0_i\in \h^0_i$, has only
zero solution. Then, $\g^*$ is the subalgebra in 
$(\n_++\g_1+\h^0_1)\oplus (\n_- + \g_2+\h^0_2)$ obtained by  
identification of $\g_1$ and $\g_2$ via $\tau$.
The classical r-matrix for $\g$ is recovered from
the canonical element of the pairing by projecting  $\g\op\g$ to the first
addend.

We restrict further considerations  to the case $\g=sl_{n-1}(\C)$; then
$D(\g)$ is represented in the direct sum
of two matrix algebras $\R=\Mat_{n}(\C)$. The cyclic inner
product is induced by the functional $t_{\R^2}=\Tr\ominus \Tr$, the 
difference of the corresponding traces.
Denote ${\Cal A}(\l)$ the associative envelope of a Lie
algebra $\l$ in $\R^2$.
Algebra ${\Cal A}(\g_1)$ is isomorphic to a direct
sum of matrix algebras $\Mat_{m_l}(\C)$.
Suppose that $\tau$ may be  extended to the isomorphism
$\hat\tau\colon {\Cal A}(\g_1) \to{\Cal A}(\g_2)$.
That means that the restriction of $\tau$  to every connected component 
of $\Gamma_1$ preserves  orientation induced  by that of the 
Dynkin diagram.
Let $\G$ be the set of diagonal
matrix idempotents $\eta_i$, $i=1,...,n$,
and  $\G_i= \G \cap {\Cal A}(\g_i)$; then
$\hat\tau$ is reduced to a bijection
$\G_1\to \G_2$. We impose one more condition  on the 
BD triple assuming
that for every $\eta \in \G_1$ there is
the smallest positive integer $m(\eta)$ such that $\hat\tau^{m(\eta)}(\eta)\not \in
\G_1$.
Again, this condition means that the subspace spanned by
$\G_1 \op \hat\tau (\G_1)$ has trivial intersection with
the subalgebra $\R\oplus \id(\R)$.
This requirement holds, for example,
if $\G_1\cap \G_2=\emptyset$ or  $\hat\tau(h_i)=h_k  \Rightarrow k>i$.

Denote $\B$ the associative subalgebra in
$({\Cal A}(\g_1)+n_+)\oplus ({\Cal A}(\g_2)+n_-)$ obtained by identification
of ${\Cal A}(\g_1)$ and ${\Cal A}(\g_2)$ via $\hat \tau$;
define  $\d_1$ and $\d_2$ as commutative subalgebras
spanned by $\G\backslash \G_1\op\{0\}$ and 
$\{0\}\op\G\backslash\G_2$ correspondingly.
Both $\d_i$ have dimension $n-m$, where $m=\card(\G_1)$;
they are orthogonal to $\B$, which is also  $\d_i$-invariant.
So the sum $\B+\d_1+\d_2$ is an associative algebra, 
and its intersection $\d$ with  ${\Cal A}(\g)$ 
is a subalgebra in   $\C\G\op\C\G$.
\proclaim{Lemma 1}
Projections $\gamma_i\colon\d \to \d_i$, $i=1,2$, have the full rank $n-m$.
\endproclaim
\vskip .03in
\demo{Proof}
We will check the statement only for $\gm_2$ in view of the
symmetry $\hat\tau \to \hat\tau^{-1}$.
Consider the sets $\G_\eta=\{\eta,\hat\tau(\eta),...,
\hat\tau^{m(\eta)}(\eta)\}$ if  $\eta\in \G_1\backslash\G_2$ and
$\G_\eta=\{\eta\}$  in case of
$\eta\in \G\backslash(\G_1\cup\G_2)$.  They do not intersect
for different $\eta$ and clearly
$\G_\eta\cap (\G\backslash\G_2)=\eta$.
The one-dimensional subspace spanned by 
$\sum_{\xi \in \G_\eta}\xi\oplus\sum_{\xi \in \G_\eta}\xi$
is evidently contained in ${\Cal A}(\g)$ which is embedded 
in $\M\oplus \M$ diagonally. It is also contained in 
$\B+\d_1+\d_2$ because its projection
to $\C\G_1\op \C\G_2$ lies in the subalgebra of $\B$
spanned by $\G_1\op \hat\tau (\G_1)$. So it is a subspace of
$\d$ and its projection ro $\d_2$ is $\C\eta$.
$\blacksquare$\enddemo
\proclaim{Corollary 1}
Algebras $\B$, and  $\d_i$, $i=1,2$, are bimodules for $\d$.
\endproclaim
\vskip .03in
\demo{Proof}
Concerning the  algebras $\d_i$, this immediately follows from Lemma 1.
By the very definition, $\d$  is a subalgebra in $\B+\d_1+\d_2$,
and the latter is  a direct sum of algebras. 
Therefore, this is also true for the algebra $\B$.
$\blacksquare$\enddemo

Introduce algebras
$$
\M=\d_2\oplus \R\oplus \R,\quad
\D=\gm_2(\d)\oplus \d,\quad
$$
$$
\N_-=\{0\}\op \R \op \id(\R),\quad
\N_+= \gm_2(\d)\op \{0\} \op \gm_2(\d) + \{0\}\op\B.
$$
The  non-degenerate cyclic inner product on
$\M$ is determined by the functional
$t_{\M}=t_{\d_2}\op t_{\R^2} =t_{\d_2}\op t_{\R}\ominus t_{\R}$,
where $t_{\d_2}$ is $t_{\R}$ restricted to ${\d_2}$.

\proclaim{Theorem 2}
Suppose $\tau$ can be extended to the automomorphism
$\hat \tau\colon {\Cal A}(\g_1)\to {\Cal A}(\g_2)$ of associative 
algebras with no 
stable points. Then, the algebras $\M$ and 
$\M_\pm=\D+\N_\pm$ form an
associative triple.
\endproclaim
\demo{Proof}
In fact, algebras $\N_\pm$ are ideals in $\M_\pm$.
That is evident for $\N_+$ and follows from 
Corollary 1 for $\N_-$. Algebras  $\N_\pm$ are 
isotropic and orthogonal to $\D$.  
It remains  to show that they are mutually dual with respect
to the inner product. That will imply 
$\D\;\cap \;(\N_- + \N_+) = \{0\}$   and  $\M=\N_- + \D + \N_+$, 
taking into account dim$(\D)\geq n-m$ (Lemma 1). We will 
compute the dual bases explicitly. The canonical element
$Q\in \N_+\oo \N_-$ includes two addends:
$$
 (0 \op e_\gm \op \theta_1(\gm) e_{\hat\tau(\gm)})
 \oo
 \sum_{k=0}^{m_2(\gm)} (0 \op f_{\hat\tau^{-k}(\gm)} \op f_{\hat\tau^{-k}(\gm)})
 ,
$$

$$
 -
 (0 \op \theta_2(\gm) f_{\hat\tau^{-1}(\gm)} \op  f_\gm)
 \oo
 \sum_{k=0}^{m_1(\gm)} (0 \op e_{\hat\tau^{k}(\gm)} \op e_{\hat\tau^{k}(\gm)})
 .
$$
Here $\gm\in\Delta^+$ and  
$e_\gm$,  $f_\gm$ are the corresponding positive and negative 
 root vectors normalized to $(e_\gm,f_{\gm})=1$ with respect
to the trace pairing. 
Functions $\theta_i$ take
two values: $\theta_i(\gm)=1$ if $\gm\in \Delta_i$ and
$\theta_i(\gm)=0$ otherwise.
The integer number $m_2(\gm)$ means the same for $\tau^{-1}$ as
$m_1(\gm)=m(\gm)$ for $\tau$. If $\gm \not \in \Delta_i$,
we assume $m_i(\gm)=0$.
There is nothing new so far in comparison with  the quasi-classical
situation, and this part of the canonical element is the same as
that of the classical r-matrix (isomorphism $\tau$ coincides 
with $\hat\tau$ on these elements). 
The distinction
appears in the sector of diagonal matrices, so we exhibit
this part of the canonical element $Q$ in  more detail:
$$
-
(\eta \op 0\op \eta)
\oo
(0 \op \eta \op \eta)
,\quad \eta \not \in \G_2\cup \G_1,
$$
$$
-
\sum_{k=0}^{m(\eta)}
(\eta \op 0\op \eta)
\oo
(0 \op \hat\tau^k(\eta) \op \hat\tau^k(\eta))
,\quad \eta \in \G_1 \backslash \G_2,
$$
$$
-
\sum_{k=0}^{m(\eta)}
(0 \op \hat\tau^{-1}(\eta)\op \eta)
\oo
(0 \op \hat\tau^k(\eta) \op \hat\tau^k(\eta))
,\quad \eta \in \G_2\cap \G_1,
$$
$$
-
(0 \op \hat\tau^{-1}(\eta)\op \eta)
\oo
(0 \op \eta \op \eta)
,\quad \eta \in \G_2\backslash\G_1.
$$
This shows that $N_+$ and  $\N_-$ are in duality.
It follows from the proof of Lemma 1 that the vectors
$$
\sum_{k,i=0}^{m(\eta)}(\eta \op  \hat\tau^{k}(\eta) \op \hat\tau^{i}(\eta))
,\quad \eta \in \G\backslash\G_2,
$$
belong to $\D$. They form an orthonormal set
and, by dimensional arguments, span entire $\D$
$\blacksquare$\enddemo

To accomplish construction of the R-matrices,
we should satisfy the Hecke condition
for some $S\in \D^{\oo2}$.
As an algebra with cyclic inner product,
$\D$ is isomorphic to $\C^{n-m}$. This problem has been solved
during the study of the standard R-matrix for $gl_n(\C)$.
Finally, to get R-matrices lying in $\Mat^{\oo2}_{n}(\C)$, let us
take the  projection  
$\pi\colon\M\to \{0\}\op \R\op \{0\}$.
The result is the sum of two terms
$$
(\pi\oo\pi)(S)=\sum_{\eta_i, \eta_j \in \G-\G_2} a^{ij}
\sum_{l=0}^{m(\eta_i)}\hat\tau^{l}(\eta_i)
\oo
\sum_{k=0}^{m(\eta_j)}\hat\tau^{k}(\eta_j),
$$
$$
(\pi\oo\pi)(Q)=
-
\sum_{\eta\in \G_1}\sum_{k=1}^{m(\eta)}
\eta
\oo
\hat\tau^k(\eta)
+
\sum_{\gm\in \Delta^+}
e_\gm
\oo
f_{\gm}
 +
\sum_{ \si,\gm \in \Delta^+,\ \si\prec\gm}
e_\gm
\wedge
f_{\si}
.
$$
Symbol $\prec$ means the partial ordering in $\Delta^+$ defined by $\tau$:
 $\si\prec\gm$ if
$\tau^k(\si)=\gm$ for some $k>0$. Numbers $a^{ik}$ satisfy
equation (6) and provide the principal solution of the Hecke condition
in the  commutative algebra $\C\G\backslash\G_2\oo\C\G\backslash\G_2$. 
Let us stress that we obtain, in this way, R-matrices
which do not tend to unit in the classical limit $q\to 1$. 
To do so,  we take $a^{ii}=-q^{-1}$ instead of $q$ for some $i$.

The deduced formula goes over into that for the standard
$gl_n(\C)$ corresponding to empty $\G_i$. Indeed, the first and the third
terms in the expression for $(\pi\oo\pi)(Q)$ vanish, and summation over
$l$ and $k$ in $(\pi\oo\pi)(S)$ is cancelled.
Another extreme possibility is the BD triple
with $\G_1=\{\eta_1,\eta_2,...,\eta_{n-1}\}$,
$\G_2=\{\eta_2,\eta_3,...,\eta_{n}\}$,
and  the isomorphism $\hat\tau\colon\eta_i \to \eta_{i+1}$.
It leads to the solution to YBE called the Cremmer-Gervais
R-matrix. In this case, the algebra $\D$ is one-dimensional, so
the  Hecke
condition is evidently fulfilled for any scalar $S=\la$. Thus we
come to the one parameter Cremmer-Gervais solution in
the form $ \la  + (\pi\oo \pi) (Q)$. This is
in agreement with \cite{H1}. Putting $\la=q/\omega$, we get
the
R-matrix of \cite{H2} for the special value of parameter $p=1$.


\subhead 6.
Remarks
\endsubhead
Analysis shows that the developed  technique  does not apply, as
it is, to the orthogonal series of simple Lie algebras.
The intuitive explanation may be the fact that rings coincide
as linear spaces with their commutator Lie algebras. Another
indication is that R-matrices for orthogonal Lie algebras do not
satisfy the Hecke condition.
So a modification of this approach for the other classical series
is an open problem. 
It has to be emphasized that associative triples are just an
algebraic construction, probably the simplest, naturally adopted for solving
the system of equations arising from decomposition of an R-matrix 
into the sum of two  solutions to YBE. Let us demonstrate
how an  extension of this scheme explains  deformation of the Yang 
R-matrix with a constant unitary  R-matrix, \cite{BFS}.
Consider a disjoint triple
$\M = \M_- + \M_+$ with $\M_pm$ isotropic and select the subspaces
$\M^c_\pm$ of ``constants''  in $\M_\pm$
annihilated by the actions $\tr$ and $\tl$:
$\M_\pm\tr \M^c_\mp = 0 = \M^c_\mp \tl \M^c_\pm$.
The sum $\M^c=\M^c_- + \M^c_+$ is a direct sum of algebras. 
In the case $\M=\Mat_n(\C)[z,\frac{1}{z}]$ considered in Section 2,
$\M^c$ is formed by constant matrices.
It can be shown that for 
every solution $S\in \M^c\oo \M^c$ to YBE  
satisfying the unitarity  condition $SS_{21} = 1\oo 1$
the sum $S+Q$ is a solution, too.

\vskip 0.5in
\end
===============================================================
Theorem (associative triples).
{\em Let $A$ and $B$ be subalgebras in $M$,
and let $D$ the orthogonal complement to $A+B$ and itself
a subalgebra.
Suppose that $A$ and $B$ are bimodules for $D$
Let $S\in D^{\ot 2}$is a solution to the YBE
Satisfying the Hecke condition 
$SS_{21}-S\si_D= \la 1$, where $\si_D$ is 
the canonical element of the restriction of
the inner product on $M$ to $D$.
Then element $R=S+Q$ is a solution to the YBE too.
}

Remark. Since $A$ is dual to $B$ w.r.t. the 
form, its restriction to the summ $A+B$ is 
non-degenerate. Hence, its restriction to 
$D$ is non-degenerate too, and $\si_D$ is
defined. 

Proof.

The YBE for the sum $R=S+Q$ reduces to 
the form
$$
S_{12}Q_{13}Q_{23}+Q_{12}S_{13}Q_{23}+Q_{12}Q_{13}S_{23}+S_{12}Q_{13}S_{23}=
$$
$$
Q_{23}Q_{13}S_{12}+Q_{23}S_{13}Q_{12}+S_{23}Q_{13}Q_{12}+S_{23}Q_{13}S_{12}
$$
because, besides $S$ and $Q$ are themselves are solutions to 
YBE,  for every $d_1, d_2\in D$ the following 
holds $(d_1\ot d_2)Q=Q(d_2\ot d_1)$. Therefore, equalities
of the type $SSQ=QSS$ are fulfilled identically and thus drop from the 
YBE for $S+Q$. Put $S=c_i\ot d^i \in D^{\ot 2}$.
and rewrite the identity above as
$$
c_i a_j \ot d^i a_k \ot b^j b^k + 
a_i c_j \ot b^i a_k \ot d^j b^k +
a_i a_j \ot b^i c_k \ot b^j d^k +
c_i a_j \ot d^i c_k \ot b^j d^k =
$$
$$
a_j c_i \ot a_k d_i \ot b^k b^j +
c_j a_i \ot a_k b_i \ot b^k d^j +
a_j a_i \ot c_k b_i \ot d^k b^j +
a_j c_i \ot c_k d_i \ot d^k b^j \>
$$

Since $A$, $B$, and $D$ covers the whole $M$, its enough
to check this identity, consecutively paring the second
tensor component with their elements.

Step 1. Paring with $x\in A$

$$
0 + a_k x c_j \ot d^j b^k + c_k x a_j \ot b^j d^k + 0 =
$$
$$
0 + c_j x a_k \ot b^k d^j + a_j x c_k \ot d^k b^j + 0 
$$

Step 2. Paring with $x\in B$

$$
c_i a_j \ot b^j x d^i + a_i c_j \ot b^j x b^i + 0 + 0 =
$$
$$
a_j c_i \ot d^i x b^j + c_j a_i \ot b^i x b^j + 0 + 0 
$$

Step 3. Paring with $x\in D$

The first and third terms on each sides turn zero.
In the last terms, perform the substitution $S\to \si_D$
of the last factor, appealing to 
the Hecke condition. The result is
$$
0 + a_i c_j \ot d^j x b^i + 0 + c_i x a_j \ot b^j x d^i =
$$
$$
0 + c_j a_i \ot b^i x d^j + 0 + a_j x c_i \ot d^i x b^j + 0 
$$
=================================================

\newpage

\subhead  Appendix. {\it Basis}
\endsubhead
\vskip 0.1in

$$
 (0 \op E_\gm \op \delta_1(\gm) E_{\tau(\gm)})
 \oo
 \sum_{k=0}^{m_2(\gm)} 0 \op E_{-\tau^{-k}(\gm)} \op E_{-\tau^{-k}(\gm)}
 ,\quad \gm\in\Delta^+
$$
$$
 -(0 \op \delta_2(\gm) E_{\tau^{-1}(\gm)} \op  E_\gm)
 \oo
 \sum_{k=0}^{m_1(\gm)} 0 \op E_{-\tau^{k}(\gm)} \op E_{-\tau^{k}(\gm)}
 ,\quad \gm\in\Delta^-
$$
Projection to the second addend gives
$$
 \sum_{\gm\in \Delta^+} E_\gm \oo E_{-\gm}
 +
 \sum_{ \si,\gm \in \Delta^+,\si\prec\gm} E_\gm \oo E_{-\si}
 -
 \sum_{\si,\gm \in \Delta^-,\si>\gm} E_\gm \oo E_{-\si}
$$
$$
 \sum_{\gm\in \Delta^+} E_\gm \oo E_{-\gm}
 +
 \sum_{ \si,\gm \in \Delta^+,\si\prec\gm} E_\gm \wedge E_{-\si}
$$

$$
-(\eta \op 0\op \eta)\oo(0 \op \eta \op \eta),
\quad \eta \not \in \G_2\cap \G_1
$$
$$
-(\eta \op 0\op \eta)\oo\sum_{k=0}^{m(\eta)}
0 \op \tau^k(\eta) \op \tau^k(\eta),
\quad \eta \in \G_1 - \G_2
$$
$$
-(0 \op \tau^{-1}(\eta)\op \eta)\oo\sum_{k=0}^{m(\eta)}
0 \op \tau^k(\eta) \op \tau^k(\eta),
\quad \eta \in \G_2\cap \G_1
$$
$$
-(0 \op \tau^{-1}(\eta)\op \eta)\oo
(0 \op \eta \op \eta),
\quad \eta \in \G_2-\G_1
$$
Projection to the second addend gives
$$
-\eta\oo\sum_{k=1}^{m(\eta)}
 \tau^k(\eta), \quad \eta \in \G_1
$$
$$
(\pi\oo\pi)S=\sum_{\eta_i, \eta_j \in \G-\G_2} a^{ij}
\sum_{l=0}^{m(\eta_i)} \tau^{l}(\eta_i)
\oo
\sum_{k=0}^{m(\eta_j)} \tau^{k}(\eta_j),
$$
$$
(\pi\oo\pi)Q=-\sum_{\eta\in \G_1}\sum_{k=1}^{m(\eta)}\eta\oo
 \tau^k(\eta)+
 \sum_{\gm\in \Delta^+} E_\gm \oo E_{-\gm}
 +
 \sum_{ \si,\gm \in \Delta^+,\si\prec\gm} E_\gm \wedge E_{-\si}
$$